\documentclass{svproc}
\usepackage{graphicx}
\usepackage{url}
\usepackage{algorithm, algorithmic}
\usepackage{amsmath,amsfonts, amssymb}   
\usepackage{booktabs} 
\newcommand{\beq}{\begin{equation}}
\newcommand{\ba}{\begin{array}}
\newcommand{\ea}{\end{array}}
\newcommand{\eeq}{\end{equation}}
\newcommand{\bdis}{\begin{displaymath}}
\newcommand{\edis}{\end{displaymath}}

\DeclareMathOperator*{\argmin}{arg\,min}

\usepackage{textcomp}
\usepackage{xcolor}
\def\BibTeX{{\rm B\kern-.05em{\sc i\kern-.025em b}\kern-.08em T\kern-.1667em\lower.7ex\hbox{E}\kern-.125emX}}

\begin{document}
\title{Total positive influence domination on weighted networks}
%
%
\author{Danica Vukadinovi\'c Greetham\inst{1}\and
Nathaniel Charlton\inst{2} \and
Anush Poghosyan\inst{3}}
\authorrunning{D. Vukadinovi\'c Greetham et al.}
%
\institute{{Knowledge Media Institute,  The Open University, UK}\\
\email{dvg27@open.ac.uk}
 \and{}
 \email{billiejoecharlton@gmail.com}
 \and {University of Bath, UK}\\
\email{A.Poghosyan.Hopes@bath.ac.uk}}
\maketitle              
\begin{abstract}
We are proposing two greedy and a new linear programming based approximation algorithm for the total positive influence dominating set problem in weighted networks. Applications of this problem in weighted settings include finding: a minimum cost set of nodes to broadcast a message in social networks, such that each node has majority of neighbours broadcasting that message;  a maximum trusted set in bitcoin network; an optimal set of hosts when running distributed apps etc. Extensive experiments on different generated and real networks  highlight advantages and potential issues for each algorithm.
\end{abstract}

\section{Introduction}\label{intro}
In  complex networks,  one often wants to control or support a dynamic process unfolding on the network. For instance, it is beneficial to have a support community  for a social network intervention to work on an individual (see e.g. \cite{Greaves:2010}), so that change affirmative messages can come from multiple sources. For this reason, an intervention designer might want to identify a support subset of the whole social network,  which will form the basis for intervention. In media and communication,  conflicting messages often travel through the network.  Identifying the minimum set of nodes that can broadcast  a particular message, so that each node hears that message  from majority of its neighbours would be useful (one can  assume that  a `broadcasting' cost is associated to each node). In distributed systems,  each node often has a cost assigned of running a distributed app. Then, the aim is to find the most cost effective subset of nodes in order to offer some resilience guarantees  to all users - for example, that each node has at least $\alpha*100\%$ of its neighbourhood  running the app.

 In graph theory, a set of nodes such that all the other, or indeed all the nodes are connected to that set is called a {\it  dominating set}. The related optimisation problem of finding a dominating set of minimum size  is NP-complete \cite{Garey:1979}. From the 1950s onward, different variants of this problem have been investigated.  Based on problems in ad-hoc communications networks,   $k$-domination was explored where each node not in the dominating set needs to have at least $k$ neighbours in the dominating set. Similarly, two versions,  \emph{positive influence dominating set} where all nodes not in the dominating set  $D$ have to have at least half of its neighbourhood in $D$ (or all nodes for the \emph{total}  version) were proposed.   A generalised version, for any percentage of neighbourhood:  $\alpha$ domination (where each node, except those in dominating set,   needs to have at least $\alpha*100$ percent of its neighbours in the dominating set) is proposed in \cite{dunbar:alphadom}. Finally, $\alpha$-rate domination is defined  in \cite{zve:upperbounds} where nodes in  the dominating set are included in the requirement, and the constraint is imposed on the closed neighbourhood (neighbours and the node itself).  Again, finding minimum cardinalities  of $\alpha$ and $\alpha$-rate dominating sets is NP-complete.
In this work,  we propose and experimentally analyse three new algorithms for the total positive influence dominating set problems on weighted networks,  thus looking into a more general problem. 

In the next section we give preliminaries and a quick overview of the relevant previous work. In Section~\ref{greedy}, we present two greedy algorithm variants based on different strategies for node selection. In Section ~\ref{algs}, an algorithm that exploits the network's community structure is proposed. We analyse the results obtained from different algorithms  ran on families of random generated graphs and real-world networks in  Section~\ref{comparison}. Finally, we discuss our results and give some pointers to future work  in Section~\ref{discussion}.

\section{Preliminaries}
Let $G(V, E)$, where $V$ is a set of nodes, and $E$ is the set of edges between them, be a simple, undirected  graph (network) with non-negative weights on the nodes, defined by a function \[w: V \to \mathbb{R^+} \cup \{0\}.\]
 We denote with $w_v$  the weight of a node $v$. 
 
Let $d_v$ denote a degree (a number of connections) for each $v \in V$. 
The neighbourhood $N(v) $ of $v$ is a set of all adjacent nodes to $v$, thus \[ N(v) =\{w | vw\in E.\}\] 
Obviously, $d_v=|N(v)|$.

A total positive influence domination set (TPIDS) $D$ is a subset of $V$ such that for each 
$v \in V$ at least  $\lceil \alpha|N(v)|\rceil$ nodes in $N(v)$ are in $D$, for  $\alpha=\frac{1}{2}$. 
We relax the last assumption and discuss the problem when $\alpha \in (0,1)$ instead.  
The weight of $D$ is $W_D=\sum_{v \in D}w_v$. We want to find the minimum weight TPIDS (MWTPIDS).

For each $v_i \in V$, $1\le i \le n$ let the variable $x_i$ has the following meaning: $x_i = 1$ if $v_i$
is contained in MWTPIDS and $x_i = 0$ otherwise. We consider the following linear programming relaxation $LP$ of an integer program $IP$ that describes MWTPIDS  (see \cite{greetham2014}):
\begin{eqnarray}\label{lpr_weighted_dom}
\min  & \quad  & \sum_{i=1}^n{w_ix_i}\\ \nonumber
\mbox{s.t.} & \quad  & \sum_{v_j\in N(v_i)}x_j\ge \lceil \alpha| N(v_i)|\rceil, \quad \forall v_i \in V\\
&\quad& 0\le x_i \le 1, \quad \forall 1\le i \le n. \nonumber
\end{eqnarray}
The $LP$ \ref{lpr_weighted_dom} is polynomial-time solvable and we can compute an optimal solution $\{\bar{x_i}\}_{1\le i \le n}.$ If we denote
with $IP_{OPT}$ an optimal solution of the corresponding integer program IP we have that
\begin{equation} \label{opt_sol}
IP_{OPT} \ge \sum_{i=1}^n w_i\bar{x_i}.
\end{equation}

\subsection{Previous work}
Due to their suitability to a wide range of applications in networks design and control, variants of domination problems have been studied thoroughly. This includes a study of corresponding computational complexities for different variants and development of exact and approximation algorithms (e.g. \cite{Klasing:2004}, \cite{amb06}, \cite{vaz18}).
The widely explored variants include the basic dominating set problem and its weighted version where weights are on nodes. The minimum weighted dominating set problem is one of the classic NP-hard optimisation problems in network theory \cite{Garey:1979}. 
A generalisation of the domination set problem on node-weighted networks, where the direct connections are replaced with shortest paths corresponding to some measure $f$ defined on the nodes of a network, was explored in
 \cite{chen:approx}. The authors have used randomised rounding to prove the approximation ratio of $O(\log\Delta')$ for their algorithm, where $\Delta'$ is the maximum cardinality of the sets of nodes that can be dominated by any single node through the defined shortest paths.
Molnar et al. \cite{molnar:2014} proposed probabilistic dominating set selection strategies for large heterogeneous non-weighted graphs and explored how the structure of networks influences performances of degree dependent probabilistic method based approximation algorithms and greedy algorithms.

 Another generalisation, the $\alpha$-domination problem, was introduced by Dunbar et al.\ in \cite{dunbar:alphadom}, where each node not in the dominating set is required to have at least $\alpha*100$ percents of neighbours in the dominating set. Similarly,  the concept of $\alpha$-rate domination \cite{zve:upperbounds} requires each node  in the network to have at least $\alpha*100$ percents of the closed neighbourhood in the dominating set. Both the  $\alpha$ and  $\alpha$-rate domination problems are proven to be NP-complete.  New upper bounds and randomised algorithms for finding the $\alpha$ and $\alpha$-rate domination sets in terms of a parameter $\alpha$ and network node degrees on undirected simple finite graphs are provided by using the probabilistic method in  \cite{zve:randomized} and \cite{zve:upperbounds}.

Wang et al.\ \cite{Wang:2012} investigated the propagation of influence in the context of social networks. They introduced new variants of domination such as the positive influence dominating set (PIDS) and total positive influence dominating set (TPIDS). Actually, the definitions of PIDS is equivalent to $\alpha$-dominating set problem for a special case when $\alpha = 1/2$.  Dinh et al.\ \cite{Dinh:2014} have generalised PIDS and TPIDS by allowing any $0<\alpha<1$,  presenting a linear time exact algorithm for trees, and approximation algorithms for minimum PIDS and TPIDS within a factor $\ln \Delta + O(1)$, where $\Delta$ is the maximum degree of the network

 A new greedy algorithm for minimum TPIDS is proposed  in \cite{dha15} and compared with two previous greedy strategies, noting that  for PIDS and TPIDS different strategies  are needed. Minimum TPIDS (as defined here, although the authors are calling it PIDS) is presented as an integer linear program in \cite{lin18}.
In  \cite{greetham2014}, the alpha-rate domination on node-weighted networks is investigated. An algorithm based on randomised rounding of linear programming formulation of the problem is given, with a proof of its approximation ratio to be $log_2 \Delta(G)$, where $\Delta$ is the maximum degree of the network.
With a slight modification, this algorithm can be applied to MWTIPDS, but for large networks, solving LP takes time \cite{lee15}. Here we contribute with three fast alternatives and compare the quality of results for different network types.
\section{Greedy algorithms}\label{greedy}
In this section, we consider two greedy techniques for solving the minimum weighted total positive influence dominating set problem. 
The Algorithm~1 below describes a generic greedy algorithm to find MWTPIDS, where we assign a cost function $g$ (defined on the next page) to all nodes according to their weight and degree, and select the nodes with minimal cost to be in the dominating set.
\begin{algorithm}
\caption{Greedy algorithms for MWTPIDS}
\begin{algorithmic}[1]
\REQUIRE A network $G$, a real number $\alpha$, $0<\alpha\le 1$
\ENSURE A low weight TPIDS  $D$ of $G$
\STATE Initialize $D=\emptyset$; \COMMENT{ Form a set $D\subseteq V(G)$}
 \FORALL{nodes $v\in G \setminus D$}
   	\STATE compute  $c_v$; assign $g(v):=w_v/c_v$
  \ENDFOR
 \WHILE { $\exists v\in V(G)$ s.t. $r>0, r:=\lceil\alpha |N(v)| \rceil - |N(v)\cap D|$} 
          \STATE  sort $g(v)$ 
        \FOR { $k<= r$} 
       		 \STATE add  $\argmin_v g(v)$ to $D$ \COMMENT{ Add smallest $r$ nodes according to $g$ }
         \ENDFOR
         \FORALL{nodes $v\in G \setminus D$}  
                 	 \STATE recompute $g(v)$  \COMMENT{ Repeat  2-4 }
           \ENDFOR
\ENDWHILE
\RETURN $D$; \COMMENT{  $D$ is a low-weight $\alpha$-rate dominating set}
\end{algorithmic}
\end{algorithm}

As expected with a greedy process, this does not necessarily yield the optimal solution.
We consider and implement two different strategies of cost computation, that determines the nodes  added to the dominating set. Those are:
\begin{description}
\item[$S1$:] $c_v=d_v-n$, where $n$ is the number of $v$'s neighbours that are in $D$ or dominated, and the cost $g(v)=\frac{w_v}{c_v}$, thus the  nodes with large degree (discounting for neighbours that are already in domination set or dominated), but smaller  weight are picked up first;
\item[$S2$:]$c_v=\sum_{u \in N(v)\setminus{D\cup S}} w_u$, $g(v)=\frac{w_v}{c_v}$ ($S$ are the dominated nodes), hence the nodes with smaller weight and large neighbourhood weight  are picked up first (again, we discount for neighbours already in $D$ or dominated).
\end{description}
 The first strategy $S1$  tries to balance minimising weight with the minimising the size of a dominating set. The second strategy $S2$ is based on reasoning that it should be beneficial to take `light' nodes with `heavy' neighbourhoods as then less heavy neighbours will be needed in the dominating set. In both cases when we calculate cost-functions, we are not considering nodes already in dominating set or dominated. 

 Since we may need to browse through all the neighbours of nodes in $V$, in total it can take $O(n^{2})$ steps to calculate domination rate for each node $v\in V(G)$. Then computing and sorting a cost function $g$ for each node can take $O(n^{2})$ steps in the worst case. This needs to be recomputed in each loop iteration, hence, in total, the set $D$ can be computed in $O(n^{3})$ steps.

\section{An algorithm using community structure - RRWC}
\label{algs}

As the range and size of a network determine the size of the linear programme that needs to be solved, we investigated the so called block separability strategy. Our problem in its linear programming form \ref{lpr_weighted_dom}, similar to the one in \cite{greetham2014}, when there is a genuine community structure in the network, can benefit from  block separability (see \cite{boy11}). The community structure here means that almost all the edges belonging to the nodes inside a community are toward the nodes in the same community and only very few are to the nodes in other communities. Therefore, the cost function can be (approximately) separated across the communities, the adjacency matrix is (approximately) block-diagonal, so we can solve LPs similar to those in \ref{lpr_weighted_dom} separately for each community (note that this can be done in parallel).   Firstly,  we split a network into communities, then we solve a linear programme for each of communities, and use randomised rounding inside communities. Finally,  we check if all the nodes are $\alpha$ dominated, and if not, we add a required number of  nodes to the final solution. We implemented this algorithm in Python using NetworkX, and the module Community\footnote{\url{http://perso.crans.org/aynaud/communities/}}
 that deploys  the Louvain method of community detection given in \cite{Blondel:2008}.

  We will denote this algorithm {\bf AlgRRWC} and it is presented below, see Algorithm~2 (RRWC stands for \emph{randomised rounding with communities}).

\begin{algorithm}[h!]
\label{algrw}
\caption{ Algorithm RRWC for MWTPIDS}
\begin{algorithmic}[1]
\REQUIRE A network $G$, a real number $\alpha$, $0<\alpha\le 1$
\ENSURE A low weight TPIDS $D$ of $G$
\STATE Initialize $D=\emptyset$;violation$=1$
\STATE Split $G$ into communities $C_1, \ldots C_k$; 
 \FORALL {$C_i$ }
 \WHILE {no-of-runs $< \lceil\log_2 \Delta(G)\rceil$ and violation$==1$}

    \STATE solve LP; $x=$lp.result;
        \FORALL {$x_i$}
        \STATE $r=random.uniform(0, 0.5)$
        \IF {$r<x_i$}
            \STATE add $x_i$ to $D$
        \ENDIF
        \ENDFOR
        \FORALL {$x_i$}
         \IF {$|D \cap N(v_i)|< \lceil\alpha*N(v_i)\rceil$}
        \STATE violation$=1$
        \ELSE
        \STATE violation$=0$
           \ENDIF
        \ENDFOR

        no-of-run++
 \ENDWHILE
\ENDFOR
  \FORALL{node $v\in G$}
         \IF {$r=\lceil\alpha*|N(v)|\rceil- |D \cap N(v_i)|\ge 0$}
        \STATE add the first lightest $r$ neighbours not already in $D$ to $D$
           \ENDIF
        \ENDFOR

\RETURN $D$; \COMMENT{  $D$ is a low-weight $\alpha$-rate dominating set}
\end{algorithmic}
\end{algorithm}

\section{Results}
\label{comparison}
We have applied our algorithms on three real-life networks obtained from Facebook and Bitcoin Alpha. In addition, to thoroughly test advantages and potential issues  for each algorithm,  we created three types of random networks and several graph colouring benchmark random graphs. All experiments were run  on a MacBookPro, MacOS High Sierra 10.13.6  with  Intel Core i7 at 3.5GHz and 16GB of RAM.
\subsection{Real-life networks}\label{real-net}
Firstly, two publicly available Facebook networks from network repository were obtained \cite{ros15}. The first one, socfb-mich67 (fb1) is quite dense without obvious community structure, while  the other one, socfb-nips-eg (fb2) is  relatively sparse connecting different ego-nets. The largest connected component is extracted in each case and  a random integer from $1$ to $10$ is assigned  as a weight to each node.
We also used a real network from \cite{sna14},
Bitcoin-alpha network `who-trusts-whom network of people who trade using Bitcoin on a platform called Bitcoin Alpha'. This is a directed network with integer edge 
weights from $-10$  (total distrust)  to $+10$ (total trust). Pre-processing was needed, in order to have positive weights only, and to turn maximisation into minimisation problem. We added $10$ to each edge score, and then calculated the node weights $w_v$  in following way: for each node $v$, $W_v$ denoted the sum of  all  
incoming edges' weights, then \[w_v=1-\frac{W_v}{max_{v \in V}(W_v)},\] so that all weights are between $0$ and $1$.  We then converted the edges to undirected.  In this way, a minimum weight (most trusted) set of nodes is obtained, such that each node had at least $\alpha$ neighbours in that set.

\begin{table*}[h!]
\center
\caption{Real-world network statistics: $V$ denotes the number of nodes, $E$ the number of edges,  $\Delta$ the max degree, $\delta_{avg}$ the average degree, Comms the number of communities detected by Louvain algorithm and  $w_{avg}$ is the average node weight.}

\begin{tabular}{llllllllllll}
\hline\noalign{\smallskip}
network&V &E & $\Delta$ & $\delta_{avg}$ & Comms& W$_{avg}$\\
\noalign{\smallskip}
\hline
\noalign{\smallskip}

fb1&	3745&	81901&	419&	43.74&	9	&5.46\\
fb2&	2888&	2981	&769&	2.06&	8&	5.48\\
bitcoinalpha	&3775&	14120&	511	&7.48	&22&	0.98\\
\hline
\end{tabular}

\label{stats}
\end{table*}

\subsection{Random networks}

We have generated three different types of random networks, with $10$ networks of each type.
They all had a similar number of nodes and  edges and were created using methods from the NetworkX \cite{Hagberg:2008} package. 

We also used some of random graphs  from DIMACS graph colouring  benchmarks \cite{rg1} and \cite{rg2}. Weights were assigned uniformly at random from integers between  $1$ and $10$ (including the boundaries).  
\subsubsection{Random networks, ER type}

Often used as a benchmark, our first type, ER network,  Erd\"os-R\'enyi model\cite{Bollobas:2001} is obtained by choosing uniformly at random  from a family $\mathcal{G}(n, m)$ of all possible networks on $n$ nodes with $m$ edges \cite{Bollobas:2001} resulting in a small diameter, high clustering coefficient and no genuine community structure. We used \verb"dense_gnm_random_graph" method from NetworkX with parameters $n=500, m=5000$ 
to create those networks and denote them with ER.

\subsubsection{Preferential attachment - high clustering networks, PN type}

 We used  another NetworkX method \verb"powerlaw_cluster_graph" to create networks  that result in approximate power-law degree distribution and high average clustering (we used parameters $n=500$, 
  $m=10$ random edges for each node, $0.8$ for probability of triangles)\cite{Holme:2002} again without genuine community structure. These networks are denoted with \emph{PN}.

 \subsubsection{Planted $l$-partition networks, PLP type}
 Additionally, we have created networks that consisted of several interlinked modules or communities (in our case 5 communities with equal sizes of $100$).  In these networks (also called planted $l$-partition graphs \cite{Condon:2001}) nodes in the same community or subgraph are interconnected with higher probability, in our case $p_{in}=0.18$ (this value provides each community similar to other types of networks density), and nodes of different communities are connected with much smaller probability, in our case $p_{out}=0.0001$. We used \verb"random_partition_graph" NetworkX method.  This results in networks having recognisable modular or block structure - with a lot of links inside those $5$ communities and only few links between different communities.
 
 \subsubsection{Graph colouring benchmark random graphs}
 
 We have also downloaded three random networks of different density: dsjc250-5  (an ER graph, but using $\mathcal{G}(n, p)$ family of all graphs with $n$ nodes and a probability of an edge between any two nodes $p$ with $n=250$ and $p=0.5$) from \cite{joh91} again without community structure, quite dense;  r250-1 from \cite{rg1} geometric random graph formed by
randomly placing $250$ vertices in a unit square, then putting edges between any two
vertices that are within $0.1$ distance of each other also from \cite{joh91}, relatively sparse with (local neighbourhoods) community structure.;  fpsol2-i-3 from \cite{rg2} from \cite{lew94} register allocation graphs  - a conflict graph of variables, with an edge between the two  if  they are active in the same range of code, with density between the other two, and some community structure. 
The weights for all the networks listed above were created by picking uniformly a random number from
$1$ to $10$. The descriptive statistics for these networks are given in Table \ref{random-stats}.

\subsection{Comparison}
We used Gurobi\footnote{\url{https://www.gurobi.com}}, a state-of-the art commercial mathematical programming solver through its Python interface, to obtain the exact IP solutions for smaller networks.  AlgRR from \cite{greetham2014}, slightly adapted to MWTPIDS, which considers the whole adjacency matrix inside a linear program, was used for larger networks.
We can see from Table~\ref{AlgRandom} that, as expected, AlgRRWC performs well for networks with well defined community structure (PLP graphs, fpsol2-i-3, and fb2 graph). It is much faster than AlgRR (times given include the computation of communities), in some instances three orders of magnitude, providing solutions of similar order to AlgRR and sometimes outperforming it. 
Comparing two greedy algorithms  we observe that   no simple winner  between the two emerges. AlgG\_S1 performs better than AlgG\_S2 and AlgRRWC on ER and r250 graphs, while AlgG\_S2 is better on PN  and dsjc250-5 graphs. For real networks, they are competitive with AlgRRWC only on bitcoin network for higher threshold where both perform identically.  Overall AlgRRWC times and results are competitive with two greedy algorithms, and therefore, AlgRRWC is recommended as a much faster alternative to AlgRR for larger networks.
 \begin{table*}[h!]
\center
\caption{Average statistics (from a sample of  $10$) for random generated networks and three random networks from DIMACS graph colouring benchmarks;  $V$ denotes the number of nodes, $E$ the number of edges, $\Delta$ the max degree, $\delta_{avg}$ the average degree, Comms the number of communities detected by Louvain algorithm and  $w_{avg}$ is the average node weight.}

\begin{tabular}{lllllllllll}
\hline\noalign{\smallskip}
network&V &E & $\Delta$ & $\delta_{avg}$ & Comms& W$_{avg}$\\
\noalign{\smallskip}
\hline
\noalign{\smallskip}
ER-500-5000	&500	&5000	 &34.9 &	 20	&10.1&	5.48\\
PLC-500-10-0.8&	500&4874.3&	143.5&	19.5&	8.4	& 5.46\\
PP5-100-0.18-0.0001&	500	&4417.7&	30.5&	17.67&	5	&5.44\\

\hline\\
dsjc250-5-rw.gexf	&250&	15668&	147&	125.34&	6	&5.53\\
r250-1-rw.gexf	&250&	867	&13&	6.94	&14	&5.34\\
fpsol2-i-3-rw.gexf&	363&	8688&	346	&47.87&	4&	5.55\\

\hline
\end{tabular}
\label{random-stats}
\end{table*}

 \begin{table*}[h!]

\caption{Alpha-rate domination sets' average sizes (\#), average weights (W) and average running times (T) for AlgG\_S1, AlgG\_S2, AlgRRWC and AlgRR for three different types of networks (with a sample of 10 networks for each type) and three random networks from DIMACS graph-colouring benchmark set.  In bold is given the best result of the three new algorithms, and a star denotes when AlgRRWC is better than AlgRR. For smaller networks Gurobi can be treated as a ground truth.}
\center
\scalebox{0.9}{
\begin{tabular}{llllllllllllll}
&&&AlgG\_S1& && AlgG\_S2 & & &AlgRRWC & & &AlgRR&\\

\noalign{\smallskip}
\hline\noalign{\smallskip}

network& $\alpha$& \#& W  &Time(s)&  \#& W  &Time(s)&\#& W  &Time(s)& \#& W  &Time(s)\\
\noalign{\smallskip}
\hline
\noalign{\smallskip}
\small
ER-500-5000&0.25&293.7&	 \bf{1095.3}&	0.50  &  300.4 & 	1142.5 &	0.81&	269.3&	1227.9&	0.49&	224& 742&	22.10\\
ER-500-5000&0.50&401.1	&1894.8&	0.36&	400.2&	\bf{1893}&	0.50&		386.9& 1968.9&	0.47&  	339.2&1467.7&	20.56\\
ER-500-5000&0.75&472.3	& \bf{2467.2}&	0.26&	472.3& 2467.9&	0.31&	471.3& 2543.9 &	0.40&	458.7&	2386&	22.61\\
PN-500 &0.25&214.3	& 698.1&	0.41	& 201.4	&\bf{689.1}&	0.46&	234.3& 1261.8&	0.56& 	137.7&	428.8&	18.49\\
PN-500&0.50&364.5&	1615.3&	0.39&	343.7&	\bf{1527.1}&	0.34&     	353.7& 1877.7&	0.60&	265.2&	 1047&	17.16\\
PN-500&0.75&471.5&	2450.3&	0.24&471.1&\bf{2444}&	0.20&    	464.5&	2552.4&	0.56&	419.1&  2093.1&	15.35\\
PLP-5-100&0.25&318.9	&1282.2&	0.53 &	322.2&	1297.8&	0.89&	217.2&\bf{	723.3$^*$}&	0.40&        216.8& 	730.9	& 	26.43\\
PLP-5-100&0.50&403.9&	1902.2	&0.38&	398.7&	 1868.8&	0.62&   	 338.1&\bf{1454.5}&	0.41& 	336& 	1430.6& 	17.20\\
PLP-5-100&0.75&   477.7& 2500.3	&0.24& 477.4&	2498.5 &0.36&	455.2&\bf{	2348$^*$}&	0.39&	458.7 &	2386&	13.77	\\
\hline
&&&AlgG\_S1& && AlgG\_S2 & & &AlgRRWC & & &Gurobi&\\
\hline
dsjc250-5	&0.25&84&	203&	0.35&	83&	\bf{199}	&0.47&85&	227&	0.57& 	73&	   166&	6.41\\
dsjc250-5	&0.5	&148&	521	&0.24&146&\bf{509}	&0.33&156&	602&	0.54&	136& 459&	42.80\\
dsjc250-5	& 0.75	&208	&989	&0.24& 203&	\bf{949}& 	0.26&220&  1116 &	0.51&	196&888	& 	903.45\\
r250-1	&0.25&226&	1096&0.16&	227	&1106 &	0.15&	141& \bf{ 675}	&	0.27&	95&	305&	1.58\\
r250-1	&0.5	&230	&\bf{1136}	&0.11&239	&1226 &	0.11&	230&	1264	&	0.28&	150&  576&	1.88\\
r250-1& 0.75	&244	&\bf{1276}	&0.06&	244	&\bf{1276} &	0.06&	247&	1330 & 	0.19&	208&  	1013	& 	1.96\\
fpsol2-i-3	&0.25&110&	403&	0.15&	113&	405&	0.15&109&  	\bf{251}	&	0.46&	94& 183&	3.41 \\
fpsol2-i-3	&0.5	&197	&774&	0.12&	196&	772	&0.13&	198&	\bf{693}&	1.04&	184&   577& 	3.22\\
fpsol2-i-3	& 0.75	&287	&1333&	0.15&281&1303&	0.15	&283&	\bf{1285} &	0.41&	274&  	1198& 	3.32 \\
\hline
\noalign{\smallskip}
\end{tabular}
}
\label{AlgRandom}
\end{table*}

 \begin{table*}[h!]
\label{alg:rn}

\caption{Alpha-rate domination sets' sizes (\#), weights(W) and running times(T) for AlgG\_S1, AlgG\_S2, AlgRRWC and AlgRR for real-world networks and $\alpha=0.25, 0.5, 0.75$ respectively.}
\center
\scalebox{0.9}{
\begin{tabular}{llllllllllllll}

&&&AlgG\_S1 & && AlgG\_S2 & & &AlgRRWC & & &AlgRR&\\

\noalign{\smallskip}
\hline\noalign{\smallskip}

Graph& $\alpha$&\#& W  &Time(s)& \#& W  &Time(s)& \# &W  &Time(s)& \# &W  &Time(s)\\
\noalign{\smallskip}
\hline
\noalign{\smallskip}
\small
fb1	&0.25&3638	&19366&	43.59&3653&	19516	&64.09&2207&\bf{13681}&	70.22&1001&3327	&14294.23\\
fb1 & 0.5&3642	&19406&	34.75&	3645& 		19436&	39.93&2882&\bf{16771}&	50.52&1841&7607	& 8036.35\\
fb1 & 0.75&3725	&20236&	16.99&3738	& 20366	&16.92&3575&\bf{19977}&	40.87&2692	&13112&	7709.50\\
fb2	&0.25&1007	& 2236&	2.34&1007	&2236&	1.30&	876&	\bf{1760}	&1.49&824&	1612	& 123.14\\
fb2	&0.5	& 1691	& 5742&	1.21&1691&	5742&	0.66&1619	&\bf{5353} &	1.33&1572	& 5266 	&1218.46\\
fb2	&0.75&2221	& \bf{9601}&	0.55&	2221&	\bf{9601}&	0.26&	2292&10184&	1.29&2296 &	10222&	1393.59\\
bitcoinalpha &0.25&3195&3137.63&45.69 & 2978 &  2921.14& 27.96&2777 &  \bf{ 2733.91}&38.33& 1361&1315.90 & 8332.39\\
bitcoinalpha &0.5&3388&\bf{3330.19}&24.40& 3388& \bf{3330.19}&18.51&3681&   3626.38&51.62 &974 & 931.70 &23303.20 \\
bitcoinalpha &0.75&3772&\bf{3713.49}& 4.74 &3772 & \bf{3713.49}&3.23	&3775&3716.49&30.89 & 2252& 2199.55 &3714.60 \\
\hline
\noalign{\smallskip}
\end{tabular}
}
\label{relnet}

\end{table*}

\section{Conclusions}\label{discussion}
 Finding the set  of minimum weight in a  network such that each node has at least $\alpha$ percentage of its neighbours in that set  can be applied to different control and intervention problems  in distributed computing and social networks. Until now, up to our knowledge, the weighted version of  total positive influence domination problem  did not receive much attention. Our contributions consist of 
two greedy strategies and a novel linear programming based algorithm  for this problem.
We thoroughly test all proposed algorithms on a diverse set of random and real-life networks with different structures.

Splitting linear programming formulation over communities and patching up a solution offers much faster solution, when  compared with the whole  network linear program AlgRR, and as expected, produces sometime even better solutions when networks have relatively well defined community structure.  Two greedy strategies, one choosing nodes by sorting them according to  their ratio of degree and weight and the other of choosing 'lighter' nodes with 'heavier' neighbourhood, seem to be perform better than AlgRRWC only for ER and PN graphs of moderate size, but with inferior solutions. For real-life networks with community structure AlgRRWC is viable and much faster alternative to AlgRR. 

It would be interesting to explore if conditioning of matrices used in linear programming formulation can speed up solutions significantly  and to see how combination of matrix rank, density  and block-models structure influences the quality and time of solutions.  Also, further investigation of how skewed weight distributions affect different algorithms is needed.
 We hope that our results will be helpful for network analysts in many different applications across social and communication networks. Python code and networks will be  made available on GitHub\footnote{\url{https://github.com/dvgreetham/domset/}}.
 \bibliographystyle{splncs03}
\bibliography{sigproc}

\end{document}